\documentclass[11pt,a4paper]{article}

\usepackage{times}
\usepackage[utf8]{inputenc}
\usepackage{latexsym,amssymb,amsmath}
\usepackage[pdftex,colorlinks=true,breaklinks=true]{hyperref}
\usepackage{xspace}
\usepackage{tikz}
\usepackage[pdftex,colorlinks=true]{hyperref}
\usepackage{pdflscape}

\usepackage[hyperref,backend=bibtex,style=alphabetic,
            doi=false,url=false,isbn=false,firstinits=true]{biblatex}
\renewbibmacro{in:}{}
\DeclareFieldFormat{postnote}{#1}
\DeclareFieldFormat{multipostnote}{#1}
\newcommand{\Galg}{\mathbf{G}}
\newcommand{\Talg}{\mathbf{T}}
\newcommand{\Lalg}{\mathbf{L}}
\newcommand{\Malg}{\mathbf{M}}
\newcommand{\Halg}{\mathbf{H}}
\newcommand{\Palg}{\mathbf{P}}
\newcommand{\Ualg}{\mathbf{U}}
\newcommand{\cE}{\mathcal{E}}
\newcommand{\cN}{\mathcal{N}}

\newcommand{\F}{\mathbb{F}}

\newcommand{\Fq}{\mathbb{F}_q}
\newcommand{\overFq}{\overline{\mathbb{F}}_q}

\newcommand{\Z}{\mathbb{Z}}
\newcommand{\C}{\mathbb{C}}
\newcommand{\Q}{\mathbb{Q}}
\newcommand{\cl}{\mathcal{C}}

\newcommand{\Spin}{\operatorname{Spin}}

\newcommand{\SL}{\operatorname{SL}}
\newcommand{\SU}{\operatorname{SU}}
\newcommand{\SO}{\operatorname{SO}}
\newcommand{\Irr}{\operatorname{Irr}}

\begin{document}

\author{Frank Lübeck\thanks{The authors acknowledges support from the DFG
within the TRR 195.}}
\title{Remarks on computing Green functions}
\date{April 3, 2020}
\maketitle

\begin{abstract}
We describe the computation of generalized Green functions and 2-parameter
Green functions for finite reductive groups.
\end{abstract}

\section{Introduction}\label{sec:intro}

Let $\Galg$ be  a connected reductive algebraic group  over an algebraic
closure $\overFq$ of a finite  field in characteristic $p$, defined over
$\Fq$,  and let  $F: \Galg  \to  \Galg$ be  the corresponding  Frobenius
morphism.

We are  interested in (complex)  class functions on  unipotent conjugacy
classes of  the finite group  $G(q) :=  \Galg^F$ of $F$-fixed  points of
$\Galg$. The  aim is to  compute explicitly  the values of  all ordinary
irreducible  characters of  $G(q)$  on unipotent  classes  (which is  an
important step  to find  the full  character table  of $G(q)$).  Here we
consider one  part of  this problem  which is  the computation  of Green
functions.

Lusztig  defined certain  class functions  on unipotent  classes, called
\emph{generalized  Green functions},  in terms  of an  induction functor
for  character sheaves.  These functions  occur in  a character  formula
for  characteristic  functions   of  character  sheaves,  see~\cite[8.3,
8.5]{MR806210}.

Deligne  and Lusztig~\cite{MR393266}  and Lusztig~\cite[1.,2.]{MR419635}
also  introduced for  $F$-stable  Levi subgroups  $\Malg$  of $\Galg$  a
twisted induction $R_{\!M}^{G}$ from class  functions on $M(q)$ to those
of $G(q)$, using  $l$-adic cohomology of certain  varieties. A character
formula  for this  map  involves certain  $2$-parameter Green  functions
which  are parameterized  by  the  unipotent classes  of  $M(q)$ and  of
$G(q)$.

In~\cite{MR1043271} Lusztig  showed (under some restrictions  on $p$ and
$q$) that  on the  level of  class functions  on unipotent  classes both
types of  induction are the  same. From this  one can conclude  that the
$2$-parameter Green functions can be computed once the generalized Green
functions are known for $M(q)$ and $G(q)$.

The generalized  Green functions can be  computed in many cases.  So, in
these cases we can compute $2$-parameter Green functions and this yields
an  essential  ingredient  for  evaluating  the  character  formula  for
$R_{\!M}^{G}$.

In  section~\ref{sec:Green}  we describe  in  more  detail the  relation
between generalized Green functions and $2$-parameter Green functions.

In  section~\ref{sec:CompGenGreen} we  collect what  is known  about the
computation of  generalized Green functions.  There are two  main steps.
The first uses the known  generalized Springer correspondence to compute
the  generalized  Green  functions  as linear  combinations  of  certain
functions which have  non-zero values only on elements lying  in a fixed
unipotent  class  of the  algebraic  group  $\Galg$.  The second  is  to
determine these functions.

In section~\ref{sec:ExSpin8} we consider  the groups $G(q) = \Spin_8(q)$
as  an example.  In that  case the  known general  results do  not fully
determine the generalized  Green functions. But we show how  to find the
$2$-parameter Green  functions for some  maximal Levi subgroups  and the
missing parameters for the generalized Green functions at the same time.

\section{Computing  2-parameter Green  functions from  generalized Green
functions}\label{sec:Green}

\subsection{$2$-parameter Green functions}\label{ssec:2pargreen}

Let $\Galg$, $F$, $G(q)$ be as in the introduction. Let $\cl$ be the set
of unipotent conjugacy classes of $G(q)$. Furthermore, let $\Malg$ be an
$F$-stable Levi subgroup of $\Galg$, $M(q) := \Malg^F$ and $\cl'$ be the
unipotent classes of $M(q)$.

Lusztig~\cite[1.,2.]{MR419635} defined  a twisted induction  from $M(q)$
to  $G(q)$ which  defines  a  linear map  $R_{\!M}^{G}$  from the  class
functions on $M(q)$ to the class functions  on $G(q)$. If $f$ is a class
function on $M(q)$ then there is  a character formula which on unipotent
classes $c \in \cl$ of $G(q)$  is (we sometimes write class functions as
functions on conjugacy classes):
\begin{equation}\label{eq:CharFormula}
(R_{\!M}^{G} f)(c) = \sum_{c' \in \cl'} 
g_{M}^{G}(c, c') f(c'),
\end{equation}
where $g_{M}^{G}(c,  c')$ is the alternating  sum of the traces  of $(u,
u') \in c \times c'$ on  certain modules defined via $l$-adic cohomology
divided by $|C_{M(q)}(u')|$ for a $u' \in c'$.

We call $g_{M}^{G}$ the $2$-parameter Green functions. For $M(q) = T(q)$
a maximal  torus and so  $c' = \{1\}$  these include the  ordinary Green
functions $Q_T^G(u) = g_T^G(u,  1)$ introduced in~\cite{MR393266}. (Note
that we use a slightly different normalization compared to the functions
$Q_M^G(u,u')  =  \frac{1}{|c'|} g_M^G(c,c')$  in~\cite[12.1]{MR1118841},
and   to    the   functions  $\gamma_M^G(u,u')   =   \frac{|M(q)|}{|c'|}
g_M^G(c,c')$ in~\cite[1.7]{MR1043271}.)

In general, these functions are difficult to compute directly from their
definition.

But  in  the   special  case  that  $\Malg  =  F(\Malg)$   is  the  Levi
complement of an $F$-stable parabolic subgroup $\Palg = \Malg\Ualg$ with
unipotent  radical $\Ualg$  the module  from $l$-adic  cohomology has  a
simple structure,  namely it is the permutation module on $G(q)/U(q)$ on
which  $G(q)$ is  acting  by  left multiplication  and  $M(q)$ by  right
multiplication, see~\cite[11.1]{MR1118841}.  So, in this case  the trace
of $(g,m) \in G(q) \times M(q)$ is just the number of fixed points under
above action. It is easy to see that with each fixed point $x U(q)$ also
$x y U(q)$ is a fixed point for all $y \in C_{M(q)}(m)$. This shows that
\begin{equation}\label{eq:splitposint}
\textrm{if }\Malg\textrm{ is contained in an $F$-stable parabolic, 
then all } g_M^G(c, c') \in \Z_{\geq 0}. 
\end{equation}

In particular, for $\Malg  = \Galg$, and so $\Ualg =  \{1\}$ we get that
$g_G^G(c, c') = \delta_{c,c'}$ for all $c, c' \in \cl$.

Conjecture: the values of $g_M^G(c, c')$ are always algebraic integers.

Furthermore,  if  we  know  a  set  of  linearly  independent  unipotent
class  functions for  $M(q)$ of  cardinality $|\cl'|$  as well  as their
images under $R_M^G$ then  we can use equation~(\ref{eq:CharFormula}) to
compute all  $g_M^G(c,c')$. Such a set  of functions is provided  by the
generalized Green functions which we recall now.

\subsection{Generalized Green functions}\label{ssec:gengreen}

Let    $\Galg$,    $F$    be    as   above.    Lusztig    has    defined
in~\cite[8.3]{MR806210} a  set of  unipotent class functions  on $G(q)$,
called generalized  Green functions. These are  partitioned into subsets
which belong  to triples $(\Lalg,C,\cE)$  of a Levi subgroup  $\Lalg$ of
$\Galg$, a  unipotent class $C$ of  $\Lalg$ and a cuspidal  local system
$\cE$  on  $C$, up  to  $\Galg$-conjugacy.  Assume  that $\Lalg$  is  an
$F$-stable Levi complement of an  $F$-stable parabolic subgroup, $F(C) =
C$ and that $F^*\cE \cong \cE$.

Fix  an  isomorphism $\tau:  F^*\cE  \to  \cE$  (unique up  to  scalar).
The  Frobenius morphism  $F$ induces  an  action on  the group  $W_\Lalg
=  N_\Galg(\Lalg)/\Lalg$.  The  $G(q)$-conjugacy classes  of  $F$-stable
$\Galg$-conjugates  of $\Lalg$  are parameterized  by the  $F$-conjugacy
classes  of  $W_\Lalg$.   For  $w  \in  W_\Lalg$  let   $\Lalg_w$  be  a
corresponding  twisted Levi  subgroup,  $C_w$ the  conjugate  of $C$  in
$\Lalg_w$ and $\cE_w$ the conjugate local system. The isomorphism $\tau$
can also  be conjugate  to an isomorphism  $\tau_w: F^*\cE_w  \to \cE_w$
(see~\cite[10.6]{MR806210}).  Lusztig defines  for each  such $(\Lalg_w,
C_w, \cE_w, \tau_w)$ (up to $G(q)$-conjugacy) a unipotent class function
$Q_{\Lalg_w, C_w, \cE_w, \tau_w}^G$ on $G(q)$.

These  functions  form   a  basis  of  the  space   of  class  functions
on   unipotent  classes   of  $G(q)$,   see~\cite[25.4]{MR849848}.  More
precisely,  they fulfill  orthogonality  relations with  respect to  the
scalar  product   $(f_1,  f_2)  :=  1/|G(q)|   \sum_{u  \in  G(q)_{uni}}
f_1(u)  \overline{f_2(u)}$ for  class  functions $f_1,  f_2$ on  $G(q)$:
the   scalar  product   of  different   Green  functions   is  $0$   and
(see~\cite[9.11]{MR806210}):
\begin{equation}\label{eq:NormGreen}
(Q_{\Lalg_w, C_w, \cE_w, \tau_w}^G, Q_{\Lalg_w, C_w, \cE_w, \tau_w}^G) =
\frac{|W_\Lalg^{wF}| \cdot q^d}{|Z^0(\Lalg)^F|},
\end{equation}
where $Z^0(\Lalg)$ is the connected center of  $\Lalg$ and $d = \dim C -
\dim \Lalg + \dim Z^0(\Lalg)$.

Under  some restriction  on the  characteristic $p$  of $\Galg$  Lusztig
showed  in~\cite[1.14]{MR1043271}   that  for  sufficiently   large  $q$
(depending on the Dynkin type of $\Galg$) we have
\begin{equation}\label{eq:GreenInduction}
Q_{\Lalg_w, C_w, \cE_w, \tau_w}^G = 
                          R_L^G (Q_{\Lalg_w, C_w, \cE_w, \tau_w}^{L_w}).
\end{equation}
The  restriction on  $p$  is  no longer  needed  because  of results  by
Shoji~\cite[I 1.11,  1.12]{MR1318530} and  Lusztig~\cite{MR3024826}. The
restriction on $q$ can at least be  dropped for the case $\Lalg = \Talg$
a maximal torus, see~\cite[II 5.5]{MR1318530}, or more generally for all
$(\Lalg,C,\cE)$ which  still occur  for the group  $\Galg/Z(\Galg)$ with
connected center, see~\cite[4.]{MR1486214}.

Now let $\Malg$  be any $F$-stable Levi subgroup of  $\Galg$. Let $(L_v,
C_v,  \cE_v,\tau_v)$  be  a  tuple  as above,  but  now  $\Lalg$  is  an
$F$-stable Levi subgroup of an $F$-stable parabolic subgroup of $\Malg$,
and $v \in N_\Malg(\Lalg)/\Lalg$.

Using equation~(\ref{eq:GreenInduction}) and the transitivity of twisted
induction (see~\cite[11.5]{MR1118841}) we get
\begin{equation}\label{eq:transitiv}
R_M^G  (Q_{L_v,  C_v,  \cE_v,\tau_v}^M)  =  R_M^G  R_L^M  (Q_{L_v,  C_v,
\cE_v,\tau_v}^{L_v}) = Q_{\Lalg_v, C_v, \cE_v, \tau_v}^G.
\end{equation}

So, if we  know the generalized Green functions of  $M(q)$ and of $G(q)$
and  if we  can identify  each Green  function $Q_{\Lalg_v,  C_v, \cE_v,
\tau_v}^G$  with  one  of  $G(q)$  (that is,  we  can  identify  the  $w
\in  W_\Lalg$ such  that $\Lalg_v$  and $\Lalg_w$  are $G(q)$-conjugate,
and  we  can  compare  $\tau_v$  and  $\tau_w$),  then  we  can  compute
the   $2$-parameter  Green  functions via~(\ref{eq:CharFormula}).    The
relations~(\ref{eq:NormGreen})  make it  easy  to invert  the matrix  of
values of generalized Green functions of $M(q)$.

\section{Computing generalized Green functions}\label{sec:CompGenGreen}

It is now possible to compute generalized Green functions in many cases,
and otherwise  to compute at least  a good ''approximation''. Here  is a
sketch of the method.

As  in the  last  section let  $(\Lalg,C,\cE)$  be a  triple  of a  Levi
subgroup  of  $\Galg$,  a  unipotent  class $C$  of  $\Lalg$  and  $\cE$
be  a  cuspidal  local  system  on $C$.  We  discuss  those  generalized
Green  functions corresponding  to  this triple  (an  important case  is
$(\Talg,\{1\},\bar{\Q_l})$ which yields the  ordinary Green functions of
$G(q)$ which are needed to compute the Deligne-Lusztig characters).

\subsection{The generalized Springer correspondence}
\label{ssec:GenSpringer}

In~\cite{MR732546} Lusztig associates a  geometric object to $\Galg$ and
$(C,\cE)$ whose decomposition into simple  objects is described in terms
of  the irreducible  characters of  the relative  Weyl group  $W_\Lalg =
N_\Galg(\Lalg)/\Lalg$. The simple object associated to a character $\chi
\in  \Irr(W_\Lalg)$ is  an object  defined by  a pair  $(C',\cE')$ of  a
unipotent class  $C'$ and  a $\Galg$-invariant irreducible  local system
$\cE'$ on $C'$.

This  construction leads  to  a bijective  map,  called the  generalized
Springer correspondence,
\begin{equation}\label{eq:Springer}
\dot{\bigcup}_{(\Lalg, C, \cE)} \Irr(W_\Lalg) \to \cN_\Galg,
\end{equation}
where $\cN_\Galg$ is  the set of pairs $(C',\cE')$ of  a unipotent class
$C'$  in  $\Galg$ and  a  $\Galg$-equivariant  irreducible local  system
$\cE'$  on  $C'$ modulo  $\Galg$-conjugacy;  the  left  hand side  is  a
disjoint  union over  all  triples  $(\Lalg, C,  \cE)$  as above  modulo
$\Galg$-conjugacy.

The   generalized   Springer    correspondence   has   been   explicitly
determined    in \cite{MR732546},    \cite{MR803339},    \cite{MR803340}
and~\cite{MR3930020}    except   for    two    small   ambiguities    in
composition   factors  of   type  $E_8$   in  characteristic   $3$  (and
in~\cite[8.6]{MR1429881} is a correction of~\cite[12.3]{MR732546}).

To understand how  the Springer correspondence can be described, we give
some more details about the sets on both sides of the map.

The set $\cN_\Galg$:  This set is a disjoint union  of the subsets which
collect the  pairs $(C',\cE')$ for  one fixed unipotent class  $C'$. So,
fix $C'$ and  $u \in C'$. Let $A_\Galg(u)  = C_\Galg(u)/C_\Galg^0(u)$ be
the  component  group  of  the  centralizer  of  $u$  in  $\Galg$.  Then
the  $\Galg$-equivariant irreducible  local systems  $\cE'$ on  $C'$ can
be  parameterized by  the irreducible  representations of  $A_\Galg(u)$.
For  simple $\Galg$  of simply-connected  type parameterizations  of the
unipotent classes and descriptions of the corresponding component groups
are known; for classical groups in  terms of partitions for Jordan block
sizes and sometimes additional invariants,  and for exceptional types in
a more ad hoc manner.

Let $\Lalg$ be a Levi subgroup of $\Galg$. Up to conjugacy we can assume
that it is a standard Levi subgroup; such that the Weyl group of $\Lalg$
is a Coxeter group with set $I$  of Coxeter generators which is a subset
of Coxeter generators $S$ of the Weyl  group $W$ of $\Galg$. In the case
that  $\Lalg$ has  a  cuspidal pair  $(C, \cE)$  as  above the  relative
Weyl group  $W_\Lalg =  N_\Galg(\Lalg)/\Lalg$ is  also a  Coxeter group,
see~\cite[Thm. 9.2]{MR732546}. It has  Coxeter generators labeled by the
complement $S \setminus I$; for $j  \in S \setminus I$ the relative Weyl
group of $\Lalg$ in the Levi subgroup corresponding to $I \cup \{j\}$ is
of  order $2$  and its  generator is  one of  the Coxeter  generators of
$W_\Lalg$. So,  using the usual  labelings of irreducible  characters of
finite Coxeter groups there is a good parameterization for the left hand
side of the Springer correspondence~(\ref{eq:Springer}).

The  articles  mentioned  above  contain explicit  descriptions  of  the
Springer  correspondence for  all  simple  simply-connected groups.  The
general  case  can be  reduced  to  these  results. First,  the  natural
quotient map  $\Galg \to \Galg/Z^0(\Galg)$  is a bijection  on unipotent
elements  and   since  $Z^0(\Galg)$   is  contained  in   the  connected
centralizer  of any  element, the  component groups  do not  change. The
description of the  generalized Springer correspondence is  the same for
$\Galg$  and  $\Galg/Z^0(\Galg)$.  So  we can  assume  that  $\Galg$  is
semisimple. Then there is a  covering map $\Galg_1 \times \Galg_2 \times
\cdots \times \Galg_k \to \Galg$ where the $\Galg_i$ ($1 \leq i \leq k$)
are simple simply-connected groups and the kernel of the map is a finite
subgroup  of  the  center  of  the left  hand  side.  Let  $\tilde{u}  =
(\tilde{u_1},  \ldots, \tilde{u_k})  \in  \Galg_1  \times \cdots  \times
\Galg_k$ be a unipotent element with image $u$ in $\Galg$, and let $A_i$
be the component  group of $\tilde{u_i}$ in  $\Galg_i$. Then $\tilde{u}$
has component  group $A_1 \times  \cdots \times  A_k$. Let $Z_i$  be the
(finite) centre of $\Galg_i$, there is  a natural map $Z_i \to A_i$, let
$\bar{Z_i}$ be its  image. The covering map is a  bijection on unipotent
elements and  maps centralizers  of unipotent elements  to centralizers,
but  the component  groups can  become smaller.  If $K  \leq Z_1  \times
\cdots \times  Z_k$ is the kernel  of the covering and  $\bar{K}$ is its
image in $\bar{Z_1}  \times \cdots \times \bar{Z_k}$  then the component
group of $u \in  \Galg$ is $A_\Galg(u) = A_1 \times  \cdots \times A_k /
\bar{K}$. So, $\Irr(A_\Galg(u))$ is in  bijection with the characters in
$\Irr(A_1  \times \cdots  \times  A_k)$ which  have  $\bar{K}$ in  their
kernel. The generalized Springer correspondence  for $\Galg$ is given by
the part of the generalized Springer correspondence of $G_1\times \cdots
\times  G_k$ where  the image  is a  pair $(C',\cE')$  for which  $\cE'$
corresponds to  one of  the characters  of the  component group  as just
described. The generalized Springer correspondence of $G_1 \times \cdots
\times G_k$  is easily described  componentwise. So, in addition  to the
description of  the generalized  Springer correspondence for  all simple
simply-connected groups  $\Galg$ the maps $Z(\Galg)  \to A_\Galg(u)$ for
all unipotent elements of $\Galg$ must be known for this reduction. This
information is implicitly  contained in the description  of the Springer
correspondence (in terms of characters of the centers).

\subsection{Generalized Springer correspondence and Frobenius action}
\label{ssec:FrobSpringer}

Now we take also the Frobenius morphism of $\Galg$ into  account. We are
only interested in  the part of the  generalized Springer correspondence
which belongs  to triples $(\Lalg,C,\cE)$ whose  $\Galg$-conjugacy class
is $F$-stable. In that case we can assume that for our representative in
this class  we have  $F(\Lalg) =  \Lalg$ is  contained in  an $F$-stable
parabolic  subgroup,  $F(C)  =  C$  and for  $u  \in  C^F$  the  induced
$F$-action on $A_\Lalg(u)$ and so $\Irr(A_\Lalg(u))$ fixes the character
corresponding to $\cE$. There is an  induced action of $F$ on $W_\Lalg$,
we  only  consider  the  generalized  Springer  correspondence  for  the
characters of $W_\Lalg$ which are $F$-stable.  Then, the image of all of
these $F$-stable characters for  all $F$-stable $(\Lalg,C,\cE)$ consists
exactly of the $F$-stable pairs $(C',\cE') \in \cN_\Galg$.

Let $\iota = (C', \cE') \in \cN_\Galg$ be an $F$-stable pair. This means
$F(C')  = C'$  so  that  there is  a  $u \in  (C')^F  =  C' \cap  G(q)$,
and  that there  is  an  isomorphism $\tau:  F^*\cE'  \to \cE'$  (unique
up  to a  nonzero  scalar).  Then $F$  induces  an  automorphism on  the
component group $A_\Galg(u)$. The Lang-Steinberg theorem yields that the
$G(q)$-conjugacy classes contained in  $(C')^F$ are parameterized by the
$F$-conjugacy classes  of $A_\Galg(u)$. If  $a \in A_\Galg(u)$  we write
$C'_a$  for the  corresponding  $G(q)$-conjugacy  class (this  labelling
depends  on  the  choice  of  $u$).  Let  $\chi$  be  the  character  of
$A_\Galg(u)$ that  parameterizes $\cE'$,  then $\chi$ is  $F$-stable and
can be  extended to the  semidirect product $A_\Galg(u)  \rtimes \langle
F\rangle $  (here $F$ means  the automorphism induced  on $A_\Galg(u)$),
this class function is nonzero and  unique on the coset $A_\Galg(u)F$ up
to multiplication with an $|F|$-th root  of unity. To $\tau$ there is an
associated  class function  $Y_\iota$ on  $G(q)$ which  is zero  outside
$(C')^F$ and has value
\begin{equation}\label{eq:Yis} 
Y_\iota (C'_a) = c \cdot \chi(aF) \textrm{ for } a \in A_\Galg(u)
\end{equation}
for some  non-zero constant $c$ which  depends on $\tau$ and  the chosen
extension of  $\chi$. See~\cite[24.2]{MR849848} or~\cite[1.3]{MR2285233}
for more details.

Let  $(\Lalg,  C, \cE)$  be  an  $F$-stable  triple  as before.  Fix  an
isomorphism  $\tau:  F^*\cE  \to  \cE$   and  fix  for  each  $F$-stable
irreducible character $\chi$ of $W_\Lalg$ an extension to the coset $W_L
F$ (in~\cite[17.2]{MR825086}  Lusztig chooses  such extensions  which he
calls preferred extensions).  Let $\iota \in \cN_\Galg$ be  the image of
$\chi$ under the Springer  correspondence. Lusztig's construction of the
generalized  Green functions~\ref{ssec:gengreen}  shows  how these  data
determine a  unique function  $Y_\iota$ as in  the last  paragraph (that
is, determine the  constant  $c$),  see again~\cite[24.2]{MR849848}  and
\cite[1.3]{MR2285233}.  The Lusztig-Shoji  algorithm  allows to  compute
the  generalized  Green functions  corresponding  to  $(\Lalg, C,  \cE)$
as  linear  combinations  of  these class  functions  $Y_\iota$,  it  is
described in~\cite[24.4]{MR849848}.  The input  data for  this algorithm
are: the  generalized Springer  correspondence, the character  values of
the extensions  of the $F$-stable $\chi$  to $W_\Lalg F$, the  orders of
the $Z^0(\Lalg_w)^F$ of the twisted groups $\Lalg_w$, and the dimensions
of the  classes $C'$ in  the image  of the Springer  correspondence (the
dimensions of  classes are  given in~\cite[I.]{MR672610}.  The algorithm
does not only  yield the Green functions but also  the precise orders of
the classes $C'_a \subset G(q)$.

What  is missing  to find  the generalized  Green functions  as explicit
class functions on $G(q)$ is to  somehow fix the isomorphisms $\tau$ and
to determine the functions $Y_\iota$ corresponding  to the choices made.
This is complicated,  but was worked out in many     cases, at least for
certain triples $(\Lalg,  C, \cE)$, by  various   authors, in particular
by  Shoji,  see~\cite{MR2285233},  \cite{MR2371771} and  the  references
given there.  Recently, Geck~\cite{Ge2020}  has determined  the ordinary
Green   functions   for  larger   rank   exceptional   types  in   small
characteristic, only  type $E_8$  was missing and  this will  be handled
in~\cite{LuE8}.

So far,  the $Y_\iota$  were not yet  determined for:  $\Spin$-groups in
odd  characteristic (the  part that  does  not appear  in $\SO$),  $\SU$
for  small  primes  (depending  on  the  rank,  for  non-ordinary  Green
functions~\cite[Thm. 4.4]{MR2285233}), exceptional groups (the part that
does not appear in the adjoint type).

For  computations  we  use  our   own  implementation  of  the  Springer
correspondence and of  the reduction to simply connected  groups. For an
overview of basic computations with root data see~\cite[2.]{MR3097233}.

\section{Examples for $G(q) = \Spin_8(q)$, $q$ odd} \label{sec:ExSpin8}

In this section we consider the untwisted simple simply connected groups
$\Galg$  of type  $D_4$ in  odd characteristic.  The fixed  point groups
under the Frobenius are isomorphic to  the spin groups $G(q) = \Galg^F =
\Spin_8(q)$.  We  will also  consider  some  Levi subgroups  $\Malg$  of
$\Galg$ of semisimple rank $3$.

We  are interested  in the  generalized  Green functions  of $G(q)$  and
$M(q)$ and  the $2$-parameter Green  functions associated to  $M(q) \leq
G(q)$.

We label the Dynkin diagram of $\Galg$ as follows:
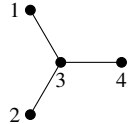
\begin{tikzpicture}[scale=0.8, baseline=-1.0]
\coordinate (3) at (0,0);
\coordinate (1) at (120:1);
\coordinate (2) at (240:1);
\coordinate (4) at (1,0);
\foreach \p in {1, 2, 3, 4} \filldraw [black] (\p) circle [radius=0.08];
\draw (1) -- (3);
\draw (2) -- (3);
\draw (3) -- (4);
\draw (1) node[anchor=east] {{\scriptsize 1}};
\draw (2) node[anchor=east] {{\scriptsize 2}};
\draw (3) node[anchor=north] {{\scriptsize 3}};
\draw (4) node[anchor=north] {{\scriptsize 4}};
\end{tikzpicture}

Corresponding  to  the  graph  automorphism  which  permutes  the  nodes
$(1,2,4)$ cyclically, there is also  an automorphism of $\Galg$ of order
$3$.

We choose  this example to contribute  to a project to  compute the full
character table for the groups  $\Spin_8(q)$. Also, we wanted to compare
the method described here  to another approach described in~\cite{RM20}.
That  alternative  method  is  more elementary  but  needs  some  tricky
computations with explicit unipotent elements.

\subsection{Generalized Springer correspondence}\label{ssec:genSprD4}

We  will  need  the   generalized  Springer  correspondence  for  simple
groups  of type  $A$  and $D$,  this  was determined  in~\cite{MR732546}
and~\cite{MR803339}.

For type  $D_4$ there are  four triples  $(\Lalg, C, \cE)$  to consider,
$\Lalg$ is a  maximal torus or a Levi subgroup   corresponding to one of
the subsets  of nodes $\{1,4\}$,  $\{2,4\}$ and $\{1,2\}$. In  the first
case  $C =  \{1\}$ and  $\cE$ is  constant. In  the other  cases $C$  is
the  regular unipotent  class  of $\Lalg$  and  for $u  \in  C$ we  have
$A_\Lalg(u)$  is cyclic  of order  $2$.  Here $\cE$  corresponds to  the
non-trivial character of $A_\Lalg(u)$.

$\Galg$  has  $12$  unipotent  conjugacy  classes  $C'$,  they  are  all
$F$-stable and  contain representatives $v  \in C' \cap G(q)$  such that
$F$ acts trivially on $A_\Galg(v)$.  For three classes $A_\Galg(v) \cong
Z_2\times Z_2$ is non-cyclic of order $4$, for seven classes $A_\Galg(v)
\cong  Z_2$  is  of  order  $2$,  and  for  the  remaining  two  classes
$A_\Galg(v)$  is  trivial.  So,  altogether  there  are  $28$  unipotent
conjugacy classes in $G(q)$.

The relative Weyl groups $W_\Lalg$ for $\Lalg$ not a torus are all of type
$B_2$ with Dynkin diagram 
\begin{tikzpicture}[scale=0.8, baseline=-1.0]
\coordinate (1) at (0,0);
\coordinate (2) at (1,0);
\foreach \p in {1, 2} \filldraw [black] (\p) circle [radius=0.08];
\draw (0,0.08) -- (1,0.08);
\draw (0,-0.08) -- (1,-0.08);
\draw (0.4,0) -- (0.6,0.2);
\draw (0.4,0) -- (0.6,-0.2);
\draw (1) node[anchor=north] {{\scriptsize 1}};
\draw (2) node[anchor=north] {{\scriptsize 2}};
\end{tikzpicture}; 
here the left node 1 corresponds in  the Dynkin diagram of $D_4$ to node
2, 1,  3, when  $\Lalg$ is  defined by  the nodes  $\{1,4\}$, $\{2,4\}$,
$\{1,2\}$,  respectively (note  that  this is  not  symmetric under  the
triality  automorphism,  the  third  case  corresponds  to  the  $\SO_8$
quotient of $\Galg$ and the other two  to the two half spin quotients of
$\Galg$ which in rank $4$ are isomorphic to $\SO_8$).

The Frobenius  induces the identity map  on all $W_\Lalg$. So  we get 13
(ordinary) Green functions and three times 5 generalized Green functions
for the cases where $\Lalg$ is of type $A_1+A_1$.

\subsection{The $Y_\iota$-functions}\label{ssec:YiD4}

The spin  group $\Galg$ has  a non-cyclic center  of order $4$.  For any
subgroup $Z =  \langle z\rangle$ of order $2$ of  this center the factor
group is $\Galg/Z \cong \SO_8(\overFq)$.

For groups  of  type  $\Halg =  \SO_N$  Shoji~\cite[Thm. 4.2]{MR2371771}
showed the  following: First he chooses for each unipotent class  $C'$ a
specific $H(q)$-class  $C'_1 \subseteq (C')^F$ (its  elements are called
split elements).  For a cuspidal pair  in $(\Lalg, C, \cE)$  one can now
normalize the isomorphism  $\tau: F^*\cE \to \cE$ by  the condition that
the characteristic function of $\tau$ has a specific positive real value
on  split elements.  Shoji shows  that in  this case  the values  of all
$Y_\iota$ (where $\iota = (C', \cE')$ is in the image of the $(\Lalg, C,
\cE)$-part  of  the  Springer  correspondence)  is  $\chi(1)$  on  split
elements  of  $C'$  where  $\chi \in  \Irr(A_\Halg(u))$  corresponds  to
$\cE'$. In  other words, when choosing  $u \in (C')^F$ as  split element
then the scalar in equation~(\ref{eq:Yis}) is always $c = 1$.

Unfortunately, a similar  description is not available  for spin groups.
But we can  deduce much from the $\SO_8$-case. If  for a unipotent class
$C'$ in  $\Galg$ with  $u \in C'$  the image of  $z$ in  $A_\Galg(u)$ is
trivial,  then the  image of  $u$ in  $\Galg/Z$ has  the same  component
group.  In this  case we  can  find a  split  element in  $(C')^F$ as  a
preimage of a split element in  $\Galg/Z$. If the component group of $u$
has order  $2$, then  only one  of the three  possible elements  $z$ has
trivial image in the component group, so a split element in $C'$ is well
defined. But there are also three  classes with component group of order
$4$ and in these cases the image of the center is equal to $A_\Galg(u)$.
If we choose two different elements $z$  as above then there is a unique
$G(q)$  class  in  $(C')^F$  which  is mapped  to  the  split  class  of
$\Galg/\langle z\rangle$ in both cases. It is not clear if this class is
also mapped to split elements  in $\Galg/\langle z\rangle$ for the third
possibility of $z$.  To conclude: We have a similar  statement about the
$Y_\iota$  for  $\Galg$  as  for  $\SO_8$, except  that  for  the  three
unipotent  classes  with  component  group  of  order  $4$  one  of  the
non-trivial $Y_\iota$ is only determined up to sign.

\subsection{Levi subgroups $\Malg$ of type $A_1+A_1+A_1$}

We  consider  the  Levi  subgroups   $\Malg$  which  correspond  to  the
subset  $\{1,2,4\}$ in  the  Dynkin diagram  of  $\Galg$. Here  $W_\Malg
=  N_\Galg(\Malg)/\Malg$  is   of  order  $2$,  so  that   we  have  two
$G(q)$-conjugacy  classes  of  $F$-stable  conjugates  of  $\Malg$.  The
tables  of generalized  Green functions  are  the same  for both  types,
we  get  them  from  the  generalized  Green  functions  for  $\SL_2(q)$
via  the simply  connected covering  of $\Malg/Z^0(\Malg)$  by a  direct
product  of three  copies  of  $\SL_2$. The  triples  $(\Lalg, C,  \cE)$
with  cuspidal pairs  are  the same  as  for $\Galg$,  so  $\Lalg$ is  a
maximal torus  or one  of the subgroups  corresponding to  the subgraphs
$\{1,4\}$,  $\{2,4\}$,  $\{1,2\}$  of  the  Dynkin  diagram.  We  remark
that  the  split element  chosen  to  normalize the  isomorphism  $\tau:
F^*\cE \to  \cE$ is  the same  in the  cases of  $\SL$ and  $\Spin$ (the
class  of  $\left(\begin{array}{cc}1&1\\0&1\end{array}\right)$  in  each
$\SL_2$-component). The Weyl  group of $\Malg$ is  elementary abelian of
order $8$,  so we get  $8$ ordinary  Green functions. The  relative Weyl
group $N_\Malg(\Lalg)/\Lalg$  is of  order $2$, so  for each  $\Lalg$ of
type $A_1+A_1$ we get $2$  generalized Green functions. Together we have
$14$ unipotent classes and generalized Green functions for each $M(q)$.

The matrix  of values of generalized  Green functions is easy  to invert
using~(\ref{eq:NormGreen}).  To find  the $2$-parameter  Green functions
with~(\ref{eq:transitiv})  we  need  to  determine  the  fusion  of  the
$F$-conjugacy classes  in $N_\Malg(\Lalg)/\Lalg$ into  the $F$-conjugacy
classes of $W_\Lalg$ in each case.

We computed the 2-parameter Green  functions first for the split version
of $\Malg$  and introduced  three indeterminates $a_{10}$,  $a_{22}$ and
$a_{27}$ for the unknown signs in three $Y_\iota$. Then we could use the
positivity and integrality  property~(\ref{eq:splitposint}) to determine
the values of these indeterminates. For example, the value
\[
g_M^G((2.2.2,4), (3221,4)) = \frac{1}{4} (a_{10}-1) (q^4+q^3)
\]
shows  that  $a_{10} =  1$  (otherwise,  the expression above  would  be
negative).

For $c, c'$ both regular there is a value
\[
\frac{1}{4}(-a_{27}+1)
\]
which  shows  that $a_{27}  =  1$  (otherwise  the  value would  not  be
integral).

And for a $c$ of type $(53)$ and $c'$ regular we find a value
\[
\frac{1}{4} (q+a_{22}-2)
\]
This shows that $a_{22}  =1$ if $q \equiv 1 \pmod{4}$  and $a_{22} = -1$
if  $q \equiv  3  \pmod{4}$. So,  the values  of  the generalized  Green
functions for $\Spin_8(q)$ depend on  the congruence class of $q$ modulo
$4$ (while everything we did so far was only depending on the congruence
modulo $2$).

To  conclude:  we  have  found   the  generalized  Green  functions  for
$\Spin_8(q)$ for  odd $q$, and  can now compute the  $2$-parameter Green
functions for the split and non-split version of $\Malg$.

\subsection{Tables}

Here are tables of $2$-parameter  Green functions. The rows are labelled
by the unipotent  classes of $M(q)$. They are described  in terms of the
Jordan block sizes in the  three $\SL_2$-components, this determines the
class in  the algebraic  group $\Malg$,  and a  number which  counts the
$M(q)$-classes of such elements; the  first corresponds  to the class of
split elements.

The columns are  labelelled by the unipotent classes of  $G(q)$. Here we
use the Jordan block sizes of images in $\SO_8$. In some cases there are
two  classes with  the  same Jordan  block sizes,  we  add a  \texttt{+}
character if the  class has representatives in the Levi  subgroup of the
subdiagram $\{1,3,4\}$ and a \texttt{-}  otherwise. Again, there is also
a  counter for  the $G(q)$-classes  in  a $\Galg$-class,  and the  first
indicates the class of split elements.

Many  values  of  the   $2$-parameter  Green  functions  factorize  with
cyclotomic polynomials as factors. In  the tables we print \texttt{P}$i$
for the $i$-th cyclotomic polynomial  evaluated at $q$ (so \texttt{P1} =
$(q-1)$, \texttt{P2} = $(q+1)$, and so on).

The first table is the case  of split $\Malg$ of type $A_1+A_1+A_1$ with
$|Z^0(\Malg)^F| = q-1$, and the second is the twisted version of $\Malg$
with $|Z^0(\Malg)^F| = q+1$.

\subsection*{Open question}

For the generalized  Green functions corresponding to  Levi subgroups of
type $A_1+A_1$ the  identification in~(\ref{eq:GreenInduction}) was only
shown for large enough $q$.

\begin{landscape}
{\scriptsize

$M(q) =A_1(q)+A_1(q)+A_1(q)+T(q-1) \leq \Spin_8(q)$, \quad $\texttt{a22} \in \{-1,1\}$ and \texttt{a22} $= q \pmod 4$

\begin{verbatim}
           | 11111111,1            221111,1 2222+,1 2222+,2 2222-,1 2222-,2 311111,1 311111,2 3221,1 3221,2 3221,3 3221,4    3311,1    3311,2
___________|_________________________________________________________________________________________________________________________________
11.11.11,1 | P2P3P4^2P6 q^4+3*q^3+3*q^2+q+1    P2P4    P2P4    P2P4    P2P4     P2P4     P2P4     P2     P2     P2     P2         1         1
 11.11.2,1 |          .               q^4P2       .       .       .       .    q^2P4    q^2P4      .      .      .      .       2*q         .
 11.2.11,1 |          .               q^4P2   q^2P4   q^2P4       .       .        .        .      .      .      .      .       2*q         .
  11.2.2,1 |          .                   .       .       . q^2P2P4       .        .        .    q^2    q^2      .      .       qP1         .
  11.2.2,2 |          .                   .       .       .       . q^2P2P4        .        .      .      .    q^2    q^2       qP1         .
 2.11.11,1 |          .               q^4P2       .       .   q^2P4   q^2P4        .        .      .      .      .      .       2*q         .
  2.11.2,1 |          .                   . q^2P2P4       .       .       .        .        .    q^2      .      .    q^2       qP1         .
  2.11.2,2 |          .                   .       . q^2P2P4       .       .        .        .      .    q^2    q^2      .       qP1         .
  2.2.11,1 |          .                   .       .       .       .       .  q^2P2P4        .    q^2      .    q^2      .       qP1         .
  2.2.11,2 |          .                   .       .       .       .       .        .  q^2P2P4      .    q^2      .    q^2       qP1         .
   2.2.2,1 |          .                   .       .       .       .       .        .        .  q^3P2      .      .      . 1/4*qP1^2 1/4*qP2^2
   2.2.2,2 |          .                   .       .       .       .       .        .        .      .      .  q^3P2      . 1/4*qP1^2 1/4*qP2^2
   2.2.2,3 |          .                   .       .       .       .       .        .        .      .      .      .  q^3P2 1/4*qP1^2 1/4*qP2^2
   2.2.2,4 |          .                   .       .       .       .       .        .        .      .  q^3P2      .      . 1/4*qP1^2 1/4*qP2^2
\end{verbatim}

\begin{verbatim}
           |    44+,1    44+,2    44-,1    44-,2   5111,1   5111,2          53,1          53,2          53,3          53,4 71,1 71,2 71,3 71,4
___________|__________________________________________________________________________________________________________________________________
11.11.11,1 |        .        .        .        .        .        .             .             .             .             .    .    .    .    .
 11.11.2,1 |        .        .        .        .        1        1             .             .             .             .    .    .    .    .
 11.2.11,1 |        1        1        .        .        .        .             .             .             .             .    .    .    .    .
  11.2.2,1 |        .        .       P2        .        .        .             1             1             .             .    .    .    .    .
  11.2.2,2 |        .        .        .       P2        .        .             .             .             1             1    .    .    .    .
 2.11.11,1 |        .        .        1        1        .        .             .             .             .             .    .    .    .    .
  2.11.2,1 |       P2        .        .        .        .        .             1             .             .             1    .    .    .    .
  2.11.2,2 |        .       P2        .        .        .        .             .             1             1             .    .    .    .    .
  2.2.11,1 |        .        .        .        .       P2        .   1/2*(a22+1)  1/2*(-a22+1)   1/2*(a22+1)  1/2*(-a22+1)    .    .    .    .
  2.2.11,2 |        .        .        .        .        .       P2  1/2*(-a22+1)   1/2*(a22+1)  1/2*(-a22+1)   1/2*(a22+1)    .    .    .    .
   2.2.2,1 | 1/2*P1P2        . 1/2*P1P2        . 1/2*P1P2        . 1/4*(q-a22-4) 1/4*(q+a22-2)   1/4*(q-a22) 1/4*(q+a22-2)    1    .    .    .
   2.2.2,2 |        . 1/2*P1P2        . 1/2*P1P2 1/2*P1P2        .   1/4*(q-a22) 1/4*(q+a22-2) 1/4*(q-a22-4) 1/4*(q+a22-2)    .    .    1    .
   2.2.2,3 | 1/2*P1P2        .        . 1/2*P1P2        . 1/2*P1P2 1/4*(q+a22-2)   1/4*(q-a22) 1/4*(q+a22-2) 1/4*(q-a22-4)    .    .    .    1
   2.2.2,4 |        . 1/2*P1P2 1/2*P1P2        .        . 1/2*P1P2 1/4*(q+a22-2) 1/4*(q-a22-4) 1/4*(q+a22-2)   1/4*(q-a22)    .    1    .    .
\end{verbatim}

\newpage

$M(q) =A_1(q)+A_1(q)+A_1(q)+T(q+1) \leq \Spin_8(q)$, \quad $\texttt{a22} \in \{-1,1\}$ and \texttt{a22} $= q \pmod 4$

\begin{verbatim}
           |  11111111,1            221111,1  2222+,1  2222+,2  2222-,1  2222-,2 311111,1 311111,2 3221,1 3221,2 3221,3 3221,4     3311,1     3311,2
___________|________________________________________________________________________________________________________________________________________
11.11.11,1 | -P1P3P4^2P6 q^4-3*q^3+3*q^2-q+1    -P1P4    -P1P4    -P1P4    -P1P4    -P1P4    -P1P4    -P1    -P1    -P1    -P1          1          1
 11.11.2,1 |           .              -q^4P1        .        .        .        .    q^2P4    q^2P4      .      .      .      .          .       -2*q
 11.2.11,1 |           .              -q^4P1    q^2P4    q^2P4        .        .        .        .      .      .      .      .          .       -2*q
  11.2.2,1 |           .                   .        .        . -q^2P1P4        .        .        .    q^2    q^2      .      .          .        qP2
  11.2.2,2 |           .                   .        .        .        . -q^2P1P4        .        .      .      .    q^2    q^2          .        qP2
 2.11.11,1 |           .              -q^4P1        .        .    q^2P4    q^2P4        .        .      .      .      .      .          .       -2*q
  2.11.2,1 |           .                   . -q^2P1P4        .        .        .        .        .    q^2      .      .    q^2          .        qP2
  2.11.2,2 |           .                   .        . -q^2P1P4        .        .        .        .      .    q^2    q^2      .          .        qP2
  2.2.11,1 |           .                   .        .        .        .        . -q^2P1P4        .    q^2      .    q^2      .          .        qP2
  2.2.11,2 |           .                   .        .        .        .        .        . -q^2P1P4      .    q^2      .    q^2          .        qP2
   2.2.2,1 |           .                   .        .        .        .        .        .        .  q^3P1      .      .      . -1/4*qP1^2 -1/4*qP2^2
   2.2.2,2 |           .                   .        .        .        .        .        .        .      .      .  q^3P1      . -1/4*qP1^2 -1/4*qP2^2
   2.2.2,3 |           .                   .        .        .        .        .        .        .      .      .      .  q^3P1 -1/4*qP1^2 -1/4*qP2^2
   2.2.2,4 |           .                   .        .        .        .        .        .        .      .  q^3P1      .      . -1/4*qP1^2 -1/4*qP2^2
\end{verbatim}

\begin{verbatim}
           |    44+,1    44+,2    44-,1    44-,2   5111,1   5111,2           53,1           53,2           53,3           53,4 71,1 71,2 71,3 71,4
___________|______________________________________________________________________________________________________________________________________
11.11.11,1 |        .        .        .        .        .        .              .              .              .              .    .    .    .    .
 11.11.2,1 |        .        .        .        .        1        1              .              .              .              .    .    .    .    .
 11.2.11,1 |        1        1        .        .        .        .              .              .              .              .    .    .    .    .
  11.2.2,1 |        .        .        .      -P1        .        .              .              .              1              1    .    .    .    .
  11.2.2,2 |        .        .      -P1        .        .        .              1              1              .              .    .    .    .    .
 2.11.11,1 |        .        .        1        1        .        .              .              .              .              .    .    .    .    .
  2.11.2,1 |        .      -P1        .        .        .        .              .              1              1              .    .    .    .    .
  2.11.2,2 |      -P1        .        .        .        .        .              1              .              .              1    .    .    .    .
  2.2.11,1 |        .        .        .        .        .      -P1   1/2*(-a22+1)    1/2*(a22+1)   1/2*(-a22+1)    1/2*(a22+1)    .    .    .    .
  2.2.11,2 |        .        .        .        .      -P1        .    1/2*(a22+1)   1/2*(-a22+1)    1/2*(a22+1)   1/2*(-a22+1)    .    .    .    .
   2.2.2,1 |        . 1/2*P1P2        . 1/2*P1P2        . 1/2*P1P2   1/4*(-q+a22) 1/4*(-q-a22-2) 1/4*(-q+a22-4) 1/4*(-q-a22-2)    1    .    .    .
   2.2.2,2 | 1/2*P1P2        . 1/2*P1P2        .        . 1/2*P1P2 1/4*(-q+a22-4) 1/4*(-q-a22-2)   1/4*(-q+a22) 1/4*(-q-a22-2)    .    .    1    .
   2.2.2,3 |        . 1/2*P1P2 1/2*P1P2        . 1/2*P1P2        . 1/4*(-q-a22-2) 1/4*(-q+a22-4) 1/4*(-q-a22-2)   1/4*(-q+a22)    .    .    .    1
   2.2.2,4 | 1/2*P1P2        .        . 1/2*P1P2 1/2*P1P2        . 1/4*(-q-a22-2)   1/4*(-q+a22) 1/4*(-q-a22-2) 1/4*(-q+a22-4)    .    1    .    .
\end{verbatim}

}
\end{landscape}
\printbibliography

\end{document}